\documentclass{amsart}
\usepackage{amsmath}
\usepackage{amssymb}
\usepackage{color}
\input xy
\xyoption{all}
\newtheorem{Lemma}{Lemma}[section]
\newtheorem{remark}[Lemma]{Remark}

\newtheorem{theorem}[Lemma]{Theorem}

\newtheorem{proposition}[Lemma]{Proposition}
\newtheorem{corollary}[Lemma]{Corollary}

\newtheorem{example}[Lemma]{Example}
\newtheorem{examples}[Lemma]{Examples}

\newenvironment{Proof}{{\sc Proof.}\ }{~$\square$\vspace{0.2truecm}}

\newcommand{\cmat}{\left(\begin{array}}
\newcommand{\fmat}{\end{array}\right)}

\newcommand{\Z}{\mathbb{Z}}

\begin{document}
   \title{Serial factorizations of right ideals}
  \author[Alberto Facchini]{Alberto Facchini}
\address{Dipartimento di Matematica, Universit\`a di Padova, 35121 Padova, Italy}
 \email{facchini@math.unipd.it}
\thanks{The first author was partially supported by Dipartimento di Matematica ``Tullio Levi-Civita'' of Universit\`a di Padova (Project BIRD163492/16 ``Categorical homological methods in the study of algebraic structures'' and Research program DOR1690814 ``Anelli e categorie di moduli''). The second author was supported by a grant from  IPM}
 \author[Zahra Nazemian]{Zahra Nazemian}
\address{School of Mathematics, Institute for Research in Fundamental Sciences (IPM), \newline P. O. Box: 19395- 5746, Tehran, Iran}
 \email{z\_nazemian@yahoo.com}
      \begin{abstract}    In a Dedekind domain $D$, every non-zero proper ideal $A$ factors as a product $A=P_1^{t_1}\cdots P_k^{t_k}$ of powers of distinct prime ideals $P_i$. For a Dedekind domain $D$, the $D$-modules $D/P_i^{t_i}$ are uniserial. We extend this property studying suitable factorizations $A=A_1\dots A_n$ of a right ideal $A$ of an arbitrary ring $R$ as a product of proper right ideals $A_1,\dots,A_n$ with all the modules $R/A_i$ uniserial modules. When such factorizations exist, they are unique up to the order of the factors. Serial factorizations turn out to have connections with the theory of $h$-local Pr\"ufer domains and that of semirigid commutative GCD domains.\end{abstract}
      
    \maketitle

\section{Introduction}

For any positive integer $n$, the factorization $n=p_1^{t_1}\cdots p_k^{t_k}$ into powers $p_i^{t_i}$ of distinct primes $p_i$ corresponds to the direct-sum decomposition $\Z/n\Z=\Z/p_1^{t_1} \Z \oplus\dots\oplus\Z/p_k^{t_k}\Z$ of the $\Z$-module $\Z/n\Z$ as a direct sum of uniserial $\Z$-modules $\Z/ \Z p_i^{t_i}$. The factorization of a non-zero integer into primes is essentially unique, and so is a direct-sum decomposition into uniserial $\Z$-modules, because uniserial modules over a commutative ring have local endomorphism rings \cite[Corollary 3.4]{Shores and Lewis}.

This fact can be generalized to the ideals of any (commutative) Dedekind domain. In a Dedekind domain $D$, every non-zero proper ideal $A$ factors as a product of powers of prime ideals in a unique way up to the order of the factors. The factorization $A=P_1^{t_1}\cdots P_k^{t_k}$ into powers $P_i^{t_i}$ of distinct prime ideals $P_i$ corresponds to the direct-sum decomposition $D/A=D/P_1^{t_1}\oplus\dots\oplus D/P_k^{t_k}$ of the $D$-module $D/A$ as the direct sum of the uniserial $D$-modules $D/P_i^{t_i}$, and these uniserial modules have local endomorphism rings.

Let $R$ be an arbitrary ring, not necessarily commutative.  In this paper, we analyze this situation considering the factorizations $A=A_1\dots A_n$ of a right ideal $A$ of $R$ as a product of proper right ideals $A_1,\dots,A_n$ with $A_iA_j=A_jA_i$ for every $i,j=1,\dots,n$, $R/A_i$ a uniserial module for every $i=1,\dots,n$ and $R/A$ canonically isomorphic to $R/A_1\oplus\dots\oplus R/A_n$. 
 We call such a factorization $A_1\dots A_n$ of $A$ a {\em serial factorization} of $A$. The endomorphism ring of a uniserial module is not local, but can have at most two maximal ideals and, correspondingly, the direct-sum decomposition of a finite direct sum of uniserial modules is not unique up to isomorphism, but depends on two permutations of two invariants, called monogeny class and epigeny class (Theorem~\ref{WKS}). Thus it is natural to expect that the corresponding factorizations $A=A_1\dots A_n$ of $A$ can be different and depend on two permutations. We show that this is not the case, and that serial decompositions, when they exist, are unique up to a unique permutation of the factors  (Theorem~\ref{vhil}).

We consider the right ideals $A$ of a ring $R$ that have a serial factorization. On the one hand, we show that if a right ideal $A$ of a ring $R$ has a serial factorization $A=A_1\dots A_n$ and $B$ is any right ideal of $R$ containing $A$, then $B$ has a serial factorization if and only if either $B\supseteq A_i$ for some index $i=1,\dots,n$ or $B$ is a two-sided ideal of $R$, and in this case we describe the serial factorization of $B$ (Theorem~\ref{Tuesday}). On the other hand, we prove that if $A,B$ are two similar right ideals of ring $R$ (that is, the right modules $R/A$ and $R/B$ are isomorphic) and $A$ has a serial factorization, then $B$ has a serial factorization and either $A=B$ or  the right $R$-module $R/A\cong R/B$ is uniserial (Theorem~\ref{3.3}).

We determine several properties of serial factorizations of right ideals, giving a number of examples and showing the analogy between the behavior of ideals with serial factorization and the behavior of  factorizations of ideals in (commutative) Dedekind rings. A commutative integral domain has the property that all its non-zero ideals have a serial factorization if and only if it is an $h$-local Pr\"ufer domain (Proposition~\ref{h-local}). We find a relation between our theory and the theory of semirigid GCD domains and factorizations of elements of a domain into rigid elements \cite{Zaf}. See Section \ref{4}.

\medskip

In this article,  $R$ is an associative ring, not necessarily commutative, with identity $1\ne 0$. Recall that a right module $M_R$ over a ring $R$ is {\em uniserial} if its lattice of submodules is linearly ordered, that is, if, for any submodules $A,B$ of $M_R$, either $A\subseteq B$ or $B\subseteq A$. A module is {\em serial} if it is a direct sum of uniserial submodules. All serial modules considered in this paper will be serial modules of finite Goldie dimension, that is, finite direct sums of uniserial submodules.

\medskip

\section{Right ideals and their serial factorizations}

Let $M$ be a right $R$-module. A finite set $\{\,N_i\mid i\in I\,\}$ of proper submodules of $M$ 
is 
{\em coindependent\/} if $N_i+ (\bigcap_{j\ne i}N_j)=M$ for
every $i\in I$, or, equivalently, if the canonical injective mapping $M/\bigcap_{i\in I}N_i\to
\oplus_{i\in I}M/N_i$ is bijective \cite[Section~2.8]{libro}. Every subfamily of a coindependent finite family of submodules is coindependent.

\begin{Lemma}\label{1} Let $A_1,\dots,A_n$ be proper right ideals of a ring $R$ such that $A_iA_j=A_jA_i$ for 
every $i,j=1,\dots,n$ and the family $\{A_1,\dots,A_n\}$ is coindependent. Then:

{\rm (1)} $A_1\dots A_n=\bigcap_{i=1}^nA_i$.

{\rm (2)} If $n \geq 2$, then each $A_i$ is a two-sided ideal. 
\end{Lemma}

\begin{Proof} (1) The proof is by induction 
on $n$. The case $n=1$ is trivial. Suppose the result true for $n-1$ proper right ideals of $R$. Let $\{A_1,\dots,A_n\}$ be a coindependent
 family of proper right ideals of $R$ such that $A_iA_j=A_jA_i$ for every $i,j=1,\dots,n$. By the inductive
  hypothesis, $A_1\dots A_{n-1}=\bigcap_{i=1}^{n-1}A_i$. Let us prove that $A_1\dots A_n=\bigcap_{i=1}^nA_i$. We have 
  that $A_1\dots A_n=A_iA_1\dots \widehat{A_i}\dots A_n\subseteq A_i$ for every $i$, so
   that $A_1\dots A_n\subseteq\bigcap_{i=1}^nA_i$. Conversely, 
   $$\begin{array}{l}\bigcap_{i=1}^nA_i=\left(\bigcap_{i=1}^nA_i\right)R=\left(\bigcap_{i=1}^nA_i\right)\left(A_n+\bigcap_{i=1}^{n-1}A_i\right)= \\ \qquad =
\left(\bigcap_{i=1}^nA_i\right)A_n+\left(\bigcap_{i=1}^nA_i\right)\left(\bigcap_{i=1}^{n-1}A_i\right)\subseteq \\ \qquad\subseteq
\left(\bigcap_{i=1}^{n-1}A_i\right)A_n+A_n\left(A_1\dots A_{n-1}\right)=A_1\dots A_n.\end{array}$$

(2) Assume that $n \geq 2$. Let $ 1 \leq i,  j  \leq n$ be such that $i \neq j$. 
Since $R = A_i + A_j$ and $A_i$ and $A_j$ commute, 
it follows that $R A_i = (A_i + A_j) A_i = A_i R = A_i$. Thus $A_i$ is two-sided.   
\end{Proof}

\begin{Lemma}\label{2.2} {\rm (1)} Let $A_1,\dots,A_n$ be $n\ge 0$ proper two-sided ideals of a ring $R$. Then  the family $\{A_1,\dots,A_n\}$ is a coindependent family if and only if $A_i+A_j=R$ for every $i\ne j$. 

{\rm (2)} If $A_1, A_2$ are two-sided ideals of $R$ and $A_1+A_2=R$, then $A_1\cap A_2=A_1A_2+A_2A_1$. 
\end{Lemma}

\begin{Proof}  (1) Suppose $A_i+A_j=R$ for every $i\ne j$. For $n=0$ and $n=1$ there is nothing to prove. Assume $n\ge2$. Then $R=(A_i+A_1)\dots\widehat{(A_i+A_i)} \dots(A_i+A_n)\subseteq A_i+(A_1\dots\widehat{A_i}\dots A_n)\subseteq A_i+\bigcap_{j\ne i}A_j$. Thus $A_i+\bigcap_{j\ne i}A_j=R$ and the family is coindependent.
The inverse implication is clear.

(2) Clearly, $A_1A_2+A_2A_1\subseteq A_1\cap A_2$. Conversely, $A_1\cap A_2=R(A_1\cap A_2)=(A_1+A_2)(A_1\cap A_2)\subseteq A_1(A_1\cap A_2)+A_2(A_1\cap A_2)\subseteq A_1A_2+A_2A_1$.\end{Proof}

\begin{proposition}\label{2} Let $R$ be a ring and $A$ be a right ideal of $R$. The following conditions are equivalent for an $n$-tuple $(A_1,\dots,A_n)$ of proper right ideals of $R$ with $A_iA_j=A_jA_i$ for every $i,j=1,\dots,n$:

{\rm (1)} $A=A_1\dots A_n$, and the family $\{A_1,\dots,A_n\}$ is a coindependent family of right ideals.

{\rm (2)} The position $r+A\mapsto (r+A_1,\dots,r+A_n)$, $r\in R$, defines a right $R$-module isomorphism $R/A\to R/A_1\oplus\dots\oplus R/A_n$.\end{proposition}

\begin{Proof} (1)${}\Rightarrow{}$(2) If the family is coindependent, the canonical mapping $R/\bigcap_{i=1}^nA_i\to
\oplus_{i=1}^nR/A_i$ is bijective. Moreover, $A=A_1\dots A_n=\bigcap_{i=1}^nA_i$ by Lemma~\ref{1}(1).

(2)${}\Rightarrow{}$(1) The position $r+A\mapsto (r+A_1,\dots,r+A_n)$, $r\in R$, defines a right $R$-module homomorphism $R/A\to R/A_1\oplus\dots\oplus R/A_n$ if and only if $A\subseteq A_i$ for every $i$. Moreover, this mapping is injective if and only if $A=\bigcap_{i=1}^nA_i$, and is an isomorphism if and only if the family $\{A_1,\dots,A_n\}$ is coindependent. We conclude by Lemma~\ref{1}(1).\end{Proof}

\begin{remark} \label{twosided}
{\rm If the equivalent conditions of Proposition \ref{2} hold for $n \geq 2$, then 
each $A_i$ is two-sided (Lemma~\ref{1}), and thus the isomorphism map in (2) is in fact a ring isomorphism. }
\end{remark}

We will now consider the right
 ideals $A$ of $R$ that have a factorization $A=A_1\dots A_n$ with $\{A_1,\dots,A_n\}$ a coindependent 
 family of proper right ideals of $R$, $A_iA_j=A_jA_i$ for every $i,j=1,\dots,n$ and $R/A_i$ a uniserial module for every $i=1,\dots,n$. 
 We call such a factorization $A_1\dots A_n$ of $A$ a {\em serial factorization} of $A$.

\begin{examples}\label{2.4} {\rm (1) Let $R$ be a commutative principal ideal domain. Then every non-zero ideal $A$ of $R$ has a serial factorization. If $A$ is generated by $a$ and $a=up_1^{t_1}\dots p_n^{t_n}$ is a factorization of $a$ with $u$ an invertible element and $p_1,\dots,p_n$ non-associate primes, then the serial factorization of $A$ is $A=A_1\dots A_n$ with $A_i=p_i^{t_i}R$. We will consider a non-commutative generalization of this example to right B\'ezout domains, that is, the integral domains in which every finitely generated right ideal is a principal right ideal, in Section~\ref{4}.

(2) More generally, let $R$ be a commutative Dedekind domain, that is, an integral domain in which every non-zero ideal factors into a product of prime ideals. Then every non-zero ideal $A$ of $R$ has a serial factorization. Namely, let $A=P_1\dots P_m$ be a factorization of $A$ into a product of prime ideals $P_i$ of $R$. Since $R$ has Krull dimension one, the non-zero prime ideals $P_i$ are maximal ideals of $R$. Thus, without loss of generality, $A=P_1^{t_1}\dots P_n^{t_n}$ with $P_1,\dots,P_n$ distinct maximal ideals of $R$. Let us show that $R/P^t$ is a uniserial module of finite composition length $t$ for every integer $t\ge0$ and every maximal ideal $P$ of $R$. The $R$-submodules of $R/P^t$ are of the form $I/P^t$ for some ideal $I$ of $R$ containing $P^t$. Then $I=Q_1\dots Q_s$ for suitable maximal ideals $Q_i$. From $I=Q_1\dots Q_s\supseteq P^t$, it follows that $P^t\subseteq Q_i$ for every $i$, so $P\subseteq Q_i$ for every $i$, hence $P= Q_i$ for every $i$. Thus $I=P^s$ for some $s$. It is now clear that the $R$-submodules $P^s/P^t$ of $R/P^t$ are linearly ordered under inclusion.}\end{examples}

Two modules $U$ and $V$ are said to have\begin{enumerate}\item {\em the same monogeny class}, denoted
$[U]_m=[V]_m$, if there exist a monomorphism $U\rightarrow V$ and
a monomorphism $V\rightarrow U$; \item {\em the same epigeny class}, denoted
$[U]_e=[V]_e$, if there exist an epimorphism $U\rightarrow V$ and
an epimorphism $V\rightarrow U$.\end{enumerate}

\begin{theorem}\label{WKS} {\rm \cite[Theorem~1.9]{TAMS}} Let $U_1$, $\dots,$ $U_n$, $V_1$, $\dots,$ $V_t$ be $n+t$ non-zero uniserial right modules over a ring $R$. Then the direct sums $U_1\oplus\dots\oplus U_n$ and $ V_1\oplus\dots\oplus V_t$ are isomorphic $R$-modules if and only if $n=t$
and there exist two permutations $\sigma$ and $\tau$ of $\{1,2,\dots,n\}$ such that $[U_i]_m=[V_{\sigma(i)}]_m$ and $[U_i]_e=[V_{\tau(i)}]_e$ for every $i=1,2,\dots, n$.\end{theorem}

\begin{proposition}\label{prev} Let $A=A_1\dots A_n$ be a serial factorization of a right ideal $A$ of a ring $R$. Then $[R/A_i]_m\ne [R/A_j]_m$  and $[R/A_i]_e\ne [R/A_j]_e$ for every $i\ne j$.\end{proposition}

\begin{Proof} The case $n = 1$ is clear. 
We first show that $[R/A_i]_m\ne [R/A_j]_m$. To see this, it suffices to prove that there is no monomorphism $R/A_i\to R/A_j$. 
Suppose the contrary. Let $\varphi\colon R/A_i\to R/A_j$ be a monomorphism. 
By Remark \ref{twosided}, $A_j$ is two-sided and so 
$A_j$ annihilates $R/A_i$. Thus  $A_j\subseteq A_i$, and therefore $R=A_i+A_j\subseteq A_i$, a contradiction.

Now we prove that $[R/A_i]_e\ne [R/A_j]_e$. It suffices to prove that there is no epimorphism $R/A_j\to R/A_i$. By contradiction, let $\psi\colon R/A_j\to R/A_i$ be an epimorphism. Since $A_j$ annihilates $R/A_j$, it also annihilates $R/A_i$, so $A_j\subseteq A_i$. 
Thus $R = A_i + A_j \subseteq A_i$, which is another contradiction. \end{Proof} 

Recall that a {\em right chain ring} is a ring $R$ with $R_R$ a uniserial right module.

\begin{theorem}\label{vhil} Let $R$ be a ring and $A$ a right ideal of $R$. Let $A=A_1\dots A_n$ and $A=B_1\dots B_m$ be two serial factorizations of $A$. Then $n=m$ and there exists a unique permutation $\sigma$ of $\{1,\dots,n\}$ such that $A_i=B_{\sigma(i)}$ for every $i=1,\dots,n$. Moreover, when $n\ge 2$,
all the right ideals $A_i$ and the right ideal $A$ are two-sided, and the factor ring $R/A$ is canonically isomorphic to the ring direct product $R/A_1\times\dots\times R/A_n$, with all the rings $R/A_i$ right chain rings.\end{theorem}

\begin{Proof} By Theorem \ref{WKS}, $n = m$.
 Suppose $n=1$. Then $m=1$, so that $A=A_1=B_1$ in this case, and the statement is trivial.

Now suppose $n\ge 2$.  By Remark \ref{twosided}, the ideals $A_i$ are two-sided. Thus $A$ is also two-sided.
By Theorem~\ref{WKS}, there is a permutation $\sigma$ of $\{1,\dots,n\}$ with $[R/A_i]_m=[R/B_{\sigma(i)}]_m$ for every $i=1,\dots,n$. Since $A_i$ and $B_{\sigma(i)}$ are two-sided ideals, the right modules $R/A_i$ and $R/B_{\sigma(i)}$ are annihilated by $A_i$ and $B_{\sigma(i)}$, respectively. From $[R/A_i]_m=[R/B_{\sigma(i)}]_m$, it follows that $A_i=B_{\sigma(i)}$.
This proves that the serial factorization of $A$ is unique up to a permutation $\sigma$ of the factors. Moreover, this permutation $\sigma$ is unique,
 because the ideals $A_i$ are coindependent and so distinct. 

Finally, since the ideals $A_i$ and $A$ are two-sided, we have that the canonical right $R$-module isomorphism $R/A\to
\oplus_{i=1}^nR/A_i$, $r+A\mapsto (r+A_1,\dots,r+A_n)$, is also a ring isomorphism, and the uniserial right $R$-modules $R/A_i$ are also right chain rings.\end{Proof}

\begin{example} {\rm The example of semisimple artinian rings shows that right ideals with a serial factorization can be rather special. Let $R$ be a semisimple artinian ring. Then all right $R$-modules are semisimple, so that, for a right ideal $A$ of $R$, $R/A$ is uniserial if and only if $A$ is a maximal right ideal of $R$. Every right ideal of $R$ is of the form $eR$ for some idempotent $e\in R$. The maximal right ideals are those of the form $(1-e)R$, where $e\in R$ is a primitive idempotent of $R$ (a {\em primitive idempotent} is an idempotent $e$ such that $eR$ is directly indecomposable). The ideals of $R$ that are finite intersections of maximal right ideals that are two-sided ideals  are those of the form $(1-e)R$, where $e\in R$ is an idempotent that is a finite sum of primitive idempotents that are central. Thus the only ideals of $R$ that have a serial factorization are those of the form $(1-e)R$, where $e\in R$ is an idempotent that is either primitive or is a finite sum of primitive idempotents that are central. 

The Wedderburn-Artin Theorem helps us in understanding how special such ideals are. Without loss of generality, we can suppose that $R:=M_{n_1}(D_1) \times \cdots \times M_{n_t}(D_t) $ for suitable integers
$t, n_1, \dots, n_t \ge 1$ and division rings $D_1, \dots, D_t$. Such a ring $R$ has infinitely many maximal right ideals (provided at least one integer $n_i$ is $>1$ and the corresponding division ring $D_i$ is infinite.) But the primitive idempotents of $R$ that are central corresponds to the indices $i$ with $n_i=1$, so that there are only finitely many two-sided ideals of $R$ with a serial factorization, possibly only the improper ideal $R$ when $n_1,\dots,n_t>1$. Cf.~Example~\ref{2.11} below.}\end{example}

\begin{example}\label{yzd} {\rm (1) Let $R_1,\dots, R_n$ be right chain rings and $R:=R_1\times\dots\times R_n$ be their direct 
product, where $n \geq 2$.
 Then every right ideal $A$ of $R$ is of the form $A:=A_1\times\dots\times A_n$, and such a right ideal $A$ is a two-sided ideal 
 of $R$ if and only if $A_i$ is a two-sided ideal of $R_i$ for every $i=1,\dots,n$. Given a proper right ideal $A:=A_1\times\dots\times A_n$ of $R$, $A$ has a serial factorization if and only if either (1) $A$ is of the form $R_1\times\dots\times R_{i-1}\times A_i\times R_{i+1}\times\dots\times R_n$ with $A_i$ a proper right ideal of $R_i$ for some $i$ or (2) $A_i$ is a two-sided ideal of $R_i$ for all $i$. In this second case the serial factorization of the right ideal $A$ is $A=\prod_{i=1}^n(R_1\times\dots\times R_{i-1}\times A_i\times\dots\times R_n)$. Here, of course, some of these $n$ factors could be the improper ideal $R$.

Notice, for instance, in the case $n=2$, that $A=A_1\times A_2$ is equal to $(A_1\times R_2)(R_1\times A_2)$ if and only if $A_2$ is a two-sided ideal of $R_2$. Also, $(A_1\times R_2)(R_1\times A_2)=(R_1\times A_2)(A_1\times R_2)$ if and only if $A_1$ is a two-sided ideal of $R_1$. 

(2) Recall that a ring is {\em right duo} if every right ideal is a two-sided ideal. Suppose that $n \geq 2$ in (1). Then 
every right ideal of $R$  has a serial factorization if and only if each right chain ring $R_i$ is right duo.
}\end{example}

\begin{example}\label{2.11} {\rm Let $S$ be any ring and $n\ge 2$ be an integer. Consider the ring $M_n(S)$ of all $n\times n$ matrices with entries in $S$. We will show that the only right ideals $I$ of $R$ with a serial factorizations are those with $R/I$ uniserial, that is, the only serial factorizations in $R$ are those with at most one factor. In fact, let $I$ be a right ideal of $R$ with a serial factorization $I=I_1\dots I_t$ with $t\ge 2$. Then $I, I_1,\dots ,I_t$ are two-sided ideals (Lemma~\ref{1}(2)). Thus they are of the form $I=M_n(J), I_1=M_n(J_1), \dots, I_t=M_n(J_t)$ for suitable two-sided ideals $J, J_1,\dots, J_t$ of $S$. All the right $R$-modules $R/M_n(J_k),$ $k=1,2,\dots,t$, must be uniserial. So, for instance the right ideals \begin{equation}\cmat{ccccc} 1&&&& \\
&0&&& \\
&&0&&\\
&&&\ddots & \\
&&&&0\fmat R +M_n(J_1)\ \ \mbox{\rm and}\ \ \cmat{ccccc} 0&&&& \\
&1&&& \\
&&0&&\\
&&&\ddots & \\
&&&&0\fmat R+M_n(J_1)\label{zeyr}\end{equation} must be comparable. But the right ideal on the left of (\ref{zeyr}) consists of all matrices with arbitrary elements of $R$ in the first row and elements of $J_1$ on the other rows, and 
the right ideal on the right of (\ref{zeyr}) consists of all matrices with arbitrary elements of $R$ in the second row and elements of $J_1$ on the other rows. These two right ideals can never be comparable. This proves that the only serial factorizations in $R$ are those with at most one factor.}\end{example}

\begin{example} {\rm Another interesting example is given by the ring $R:=T_n(k)$ of all $n\times n$ upper triangular matrices with entries in a (commutative) field $k$. Here $n\ge1$ is an integer. The ring $R$ is a serial ring, that is, both the right module $R_R$ and the left module $_RR$ are serial. For instance, if $e_{ij}$ denotes the matrix that is $1$ in the $(i,j)$ entry and $0$ in the other entries, then $R_R=e_{11}R\oplus\dots\oplus e_{nn}R$. Each right ideal $e_{ii}R$ is a uniserial $R$-module of finite composition length $n-i+1$, and its non-zero submodules are the $n-i+1$ modules $e_{ij}R$, $j=i,i+1,\dots,n$.

The ring $R$ is not a right chain ring for $n>1$, for instance because its socle is a direct sum of $n$ simple modules. The Jacobson radical $J(R)$ of $R$ consists of all strictly upper triangular matrices, so that $R/J(R)$ is the direct product of $n$ copies of $k$. It follows that $R$ has exactly $n$ maximal right ideals $M_1,\dots,M_n$, which are all left ideals as well. Here $M_i$ denotes the set of matrices in $R$ with $0$ in the $(i,i)$ entry. Notice that, for all indices $i,j=1,\dots, n$ with $i\le j$, we have that $M_i^2=M_i$ for every $i=1,\dots,n$, $M_iM_{i+1}\ne M_{i+1}M_i$ for every $i=1,\dots,n-1$, and $M_iM_j=M_jM_i=M_i\cap M_j$ for all $i,j$ with $i+1<j$.

For any two-sided ideal $A$ of $R$, we have that $A=e_{11}A\oplus\dots\oplus e_{nn}A$ as a right $R$-module, so that if $R/A$ is a uniserial right $R$-module, then there exists an index 
$ i=1 ,\dots, n$ such that $e_{11},\dots,e_{i-1\, i-1},e_{i+1\, i+1},\dots, e_{nn}\in A$. Note that $M_i$ is generated by 
all the matrices $e_{kt}$ with   $k < t$ or $k = t \neq i$.
 But 
$A$ is a two sided ideal of $R$. So, we have that, for every $k < t $, if 
$k \neq i$, then $e_{kt} = e_{kk} e_{kt} \in A$, and if $k = i$, then $t \neq i$, so 
$e_{kt} = e_{kt} e_{tt} \in A$. It follows that $A$ is equal to  $M_i$. We have thus proved that the two-sided ideals $A$ of $R$ with $R/A$ a uniserial right $R$-module are only $M_1,M_2,\dots,M_n$ and $R$.

It follows that the only two-sided ideals of $R$ that have a serial factorization are those of the form $\prod_{i\in Y}M_i$, where $Y$ is a subset of $ \{1,2,\dots,n\}$ with the property that, for every $i=1,2,\dots,n-1$, $i\in Y$ implies $i+1\not\in Y$, that is, the subsets $Y$ of $ \{1,2,\dots,n\}$ with no consecutive indices (because $M_iM_{i+1}\ne M_{i+1}M_i$). 

It can also be proved that the other right ideals $I$ of $R$ with a serial factorization, that is, those with $R/I$ a uniserial right $R$-module, are the kernels of the right $R$-module morphisms $\lambda_{e_{ii}r}\colon R\to e_{ii}R/e_{ij}R$ ($i\le j$), induced by left multiplication by an element $e_{ii}r\in e_{ii}R$. Thus they are the right ideals $(e_{ij}R:e_{ii}r):=\{\, x\in R\mid e_{ii}rx\in e_{ij}R\,\}$.
}\end{example}

\begin{example} {\rm Every right ideal of $R$ with a serial factorization is the intersection of finitely many right ideals $A_i$. Thus it is natural to ask whether the intersection of two right ideals with a serial factorization is a right ideal with a serial factorization. The answer is negative, as the following example shows. Let $S$ be a right chain ring and $I$ a right ideal of $S$ that is not two-sided. Then the right ideals $I\times S$ and $S\times I$ of $R:=S\times S$ have the trivial serial factorizations of length one, but their intersection $I\times I$ is a right ideal that is not two-sided, and $R/I\times I$ is not uniserial, so that $I\times I$ does not have a serial factorization.}\end{example}

\section{Structure of right ideals with a serial factorization} \label{3}

\begin{theorem}\label{Tuesday} Let $R$ be a ring, $A$ a right ideal of $R$ with a serial factorization $A=A_1\dots A_n$ and $B$ a right ideal of $R$ containing $A$. Then:

{\rm (1)} $B$ has a serial factorization if and only if either $B\supseteq A_i$ for some index $i=1,\dots,n$ or $B$ is a two-sided ideal of $R$.

{\rm (2)} If $B$ has a serial factorization, then the serial factorization of $B$ is $B=(B+A_1)\dots(B+A_n)$ (where we are supposed to omit the factors $B+A_i$ equal to $R$).
\end{theorem}

\begin{Proof} If $n = 1$, then (1) and (2) are clear.  Assume that $n\geq 2$. 

If $B$ has a serial factorization, then  either $R/B$ is a uniserial $R$-module or $B$ is a two-sided ideal of $R$ by Theorem~\ref{vhil}.
Since $R/B$ is isomorphic to a factor module  of $R/A_1  \oplus \dots \oplus R/A_n$,
in case $R/B$ is uniserial, it is isomorphic to a factor module of $R/A_i$ for some $ i=1,2,\dots,n$.
 Thus  $R/ B \cong R/K$ for some right ideal  $K$ of $R$ with $A_i \subseteq K$. In particular, the two-sided ideal $A_i$ annihilates $R/K$, hence annihilates $R/B$, and so $A_i \subseteq B$. This proves one of the implications in (1).
  
Let us prove the other implication. If $A_i \subseteq B$ for some $i$, then $R/B$ is  isomorphic to a factor module of 
 uniserial module $R/A_i$ and so $R/B$ 
is a uniserial module and $B = B$ is its serial factorization. Suppose
 that $B$ is a two-sided ideal of $R$. There is a canonical ring isomorphism $\varphi\colon R/A\to R/A_1\times\dots\times 
R/A_n$, $\varphi\colon r+A\mapsto (r+A_1,\dots, r+A_n)$, and all the rings $R/A_i$ are right chain rings. Since the ideals $A_i$ are 
coindependent, for every $i=1,2,\dots,n$ there exist $x_i\in A_i$ and $y_i\in\bigcap_{j\ne i}A_j$ with $1=x_i+y_i$. The inverse of the 
mapping $\varphi$ is defined by $(r_1+A_1,\dots, r_n+A_n)\mapsto \sum_{i=1}^ny_ir_i+A$. It is easily seen that in these mutually 
inverse isomorphisms, the ideal $B/A$ of $R/A$ corresponds to the ideal $B+A_1/A_1\times\dots \times B+A_n/A_n$ of 
$R/A_1\times\dots\times R/A_n$, 
the ideals $A_i/A$ of $R/A$ correspond to the ideals $R/A_1\times\dots \times A_i/A_i\times\dots\times R/A_n$, and the ideals 
$B+A_i/A$ of $R/A$ correspond to $R/A_1\times\dots \times B+A_i/A_i\times\dots\times R/A_n$. Thus the ideals $R/A_1\times\dots 
\times B+A_i/A_i\times\dots\times R/A_n$ and
 $R/A_1\times\dots \times B+A_j/A_j\times\dots\times R/A_n$ of $R/A_1\times\dots\times R/A_n$ commute, so that their corresponding 
 ideals $B+A_i/A$ and $B+A_j/A$ of $R/A$ commute. Hence $(B+A_i)(B+A_j)+A=(B+A_j)(B+A_i)+A$. 
 But $A_iA_j=A_jA_i\supseteq A$,
  so that $(B+A_i)(B+A_j)=(B+A_j)(B+A_i)$. Moreover, the
   product of the ideals $R/A_1\times\dots \times B+A_i/A_i\times\dots\times R/A_n$ of $R/A_1\times\dots\times R/A_n$ is 
$B+A_1/A_1\times\dots \times B+A_n/A_n$.
Correspondingly, in the ring $R/A$, we find that the product of the ideals $B+A_i/A$ of $R/A$ is $B/A$. Thus 
$(B+A_1)\dots(B+A_n)+A=B$. But $(B+A_1)\dots(B+A_n)\supseteq A_1\dots A_n=A$, so that $(B+A_1)\dots(B+A_n)=B$. Finally, 
since the ideals $A_1,\dots, A_n$ of $R$ are coindependent, the ideals $B+A_1,\dots, B+A_n$ are coindependent.
\end{Proof}

As  a corollary, we get that:

\begin{corollary} Let $R$ be a ring and $A,B$ be two right ideals of $R$. Let $A=A_1\dots A_n$ and $B=B_1\dots B_m$ be two serial factorizations of $A,B$ respectively. Suppose $A\subseteq B$. Then there exists an injective mapping $\sigma\colon\{1,\dots,m\}\to \{1,\dots,n\}$ such that $A_{\sigma(j)}\subseteq B_j$ for every $j=1,\dots,m$.\end{corollary}

\begin{Proof} Follows from Theorems \ref{vhil} and  \ref{Tuesday}(2).
\end{Proof}
 
 Consider the trivial case of a right chain ring $R$. Every right ideal $A$ of $R$ has the trivial unique serial factorization $A=A$ of length $1$. But a right ideal $A$ is not necessarily two-sided. Thus if in a ring $R$ all non-zero right ideals have a serial factorization, then for every right ideal $A$ of $R$, either $R/A$ is a uniserial right $R$-module or $A$ is a two-sided ideal.
As a consequence of Theorem~\ref{vhil}, we see that if $R$ is
any ring, the zero ideal of $R$ has a serial factorization if and only if 
 $R$ is a finite direct product of right chain rings.  More generally, combining this remark, Example~\ref{yzd}(2) and Theorem~\ref{Tuesday}, we get:
  
 \begin{proposition} If $R$ is a ring, every right ideal of $R$ has a serial factorization if and only if either $R$ is a right chain ring or $R$ is a finite direct product of right duo right chain rings.  
 \end{proposition}

 Recall that two right ideals $B,C$ of $R$ are {\em similar} if the right modules $R/B$ and $R/C$ are isomorphic \cite[p.~76]{Cohn}.

 \begin{theorem}\label{3.3} Let $R$ be a ring, and $A,B$ be two similar right ideals of $R$. Suppose that $A$ has a serial factorization. Then:
 
 {\rm (1)} $B$ has a serial factorization.
 
 {\rm (2)} Either $A=B$ or  the right $R$-module $R/A\cong R/B$ is uniserial.\end{theorem}
 
 \begin{Proof}Step 1: {\em If $R$ is a right chain ring, $B$ is a right ideal of $R$ and $B$ is similar to the zero ideal, then $B=0$.}
 
 Suppose $B$ similar to the zero ideal, that is, $R/B$ and $R_R$ isomorphic. Then $B$ is necessarily a proper ideal of $R$. Let $\varphi\colon R_R\to R/B$ be a right $R$-module isomorphism. Set $\varphi(1)=c+B$ for some $c\in R$. Clearly,  $cR+B=R$ and the right annihilator of the element $c+B$ of $R/B$ is zero, that is, for every $r\in R$, $cr\in B$ implies $r=0$. As $cR+B=R$ and $R$ is a local ring, we have that $cR=R$, so $c$ is right invertible, hence invertible. We must prove that $B=0$. If $b\in B$, then $c^{-1}b\in R$ and $c(c^{-1}b)=b\in B$. Thus $c^{-1}b=0$, so $b=0$, as desired.
 
 \medskip
 
Step 2: {\em If $R=R_1\times\dots\times R_n$ is a finite direct product of right chain rings $R_i$, $B$ is a right ideal of $R$ and $B$ is similar to the zero ideal, then $B=0$.}
 
Every right ideal of $R$ is of the form $B=B_1\times\dots\times B_n$ for suitable right ideals $B_i$ of $R_i$. Moreover, $R/B\cong R$ as a right $R$-module implies $R_i/B_i\cong R_i$ as right $R_i$-modules for every $i=1,2,\dots,n$.

\medskip

Step 3:  {\em Proof of {\rm (1)}.}

Suppose, on the contrary, that (1) does not hold. Thus there exists a ring $R$ with two similar right ideals $A$ and $B$, and $A$ has a serial factorization, but $B$ has not. Since $R/A\cong R/B$, the right ideal $A$ has a serial factorization and $B$ has not, it is not possible that $R/A\cong R/B$ is a uniserial $R$-module. Thus $A$ is necessarily a two-sided ideal. The modules $R/A\cong R/B$ must have the same annihilators, so that $A\subseteq B$. Now pass to the factor ring $R/A$, getting that $R/A$ is a finite direct product of right chain rings, in which the zero ideal has a serial factorization. The ideal $B/A$ of $R/A$ is non-zero (otherwise $B=A$ would have a serial factorization), and $B/A$ is similar to the zero ideal of $R/A$. This contradicts what we have proved in Step 2.

\medskip

{\em Step 4: Proof of {\rm (2)}.}

Let $R$ be an arbitrary ring, and $A,B$ be two similar right ideals of $R$. Suppose that $A$ has a serial factorization.  Then either $R/A\cong R/B$ is uniserial or $R/A\cong R/B$ is not uniserial. If $R/A\cong R/B$ is uniserial, there is nothing to prove. If $R/A\cong R/B$ is not uniserial, then both $A$ and $B$ have a serial factorization by (1), so that both $A$ and $B$ are two-sided. But then $R/A$ and $R/B$ have the same annihilators, that is, $A=B$.
 \end{Proof}
 
 \begin{proposition}\label{3.4} The following conditions are equivalent for a ring $R$.
 
 {\em (1)} Every non-zero right ideal of $R$ has a serial factorization.
 
  {\em (2)} The right ideal $rR$ has a serial factorization for every non-zero element $r\in R$.\end{proposition}
  
  \begin{Proof} It is clear that (1) implies (2). Conversely, suppose that (2) holds. Let $A$ be a non-zero right ideal of $R$. If $A$ contains 
  an element $r$ with $R/rR$ uniserial, then $R/A$ is a uniserial right $R$-module, so that $A$ has the trivial serial factorization. If $A$ 
  does not contain elements $r$ with $R/rR$ uniserial, then $rR$ is a two-sided ideal for every $r\in A$, so that $A$ is two-sided. Fix a non-zero element $r\in A$. Then $A$ has a serial factorization by Theorem~\ref{Tuesday}(1).\end{Proof}

As examples of rings satisfiying the equivalent conditions of this proposition, we have commutative principal ideal rings, commutative Dedekind domains, right chain rings and finite direct products of right duo right chain rings (for instance, finite direct products of commutative valuation rings). We will see in Proposition~\ref{h-local} that a commutative integral domain satisfies the equivalent conditions of Proposition~\ref{3.4} if and only if it is an $h$-local Pr\"ufer domain.

\medskip

Now we are going to examine some consequences of Theorem~\ref{vhil}. 
Recall that a right module $M_R$ is {\em B\'ezout} if every finitely generated submodule of $M_R$ is cyclic. A ring $R$ is {\em right B\'ezout} if the right module $R_R$ is B\'ezout.

\begin{proposition}\label{3.5} Let $R$ be a ring and $A$ a right ideal of $R$ with a serial factorization. Then $R/A$ is a B\'ezout right $R$-module.\end{proposition}

\begin{Proof} Let $A=A_1\dots A_n$ be a serial factorization of $A$. If $R/A$ is a uniserial module, then $R/A$ is B\'ezout trivially. If $R/A$ is not uniserial, then all the right ideals $A_i$ and the right ideal $A$ are two-sided, and the factor ring $R/A$ is canonically isomorphic to the direct product $R/A_1\times\dots\times R/A_n$, with all the rings $R/A_i$ right chain rings (Theorem~\ref{vhil}). Right chain rings are B\'ezout right modules over themselves, and a direct product of right B\'ezout rings is a right B\'ezout ring. Thus $R/A$ is a B\'ezout right module over itself, hence a B\'ezout right $R$-module.\end{Proof}

\bigskip

Recall that a {\em complete set of indecomposable central idempotents} of a ring $R$ is a finite set $\{e_1,\dots, e_n\}$ of pairwise orthogonal central idempotents of $R$ with $1=e_1+\dots+ e_n$ and all the idempotents $e_i$ centrally primitive (a {\em centrally primitive idempotent} is a central idempotent that cannot be written as the sum of two non-zero orthogonal central idempotents). Complete sets of indecomposable central idempotents of $R$ correspond to direct-product decompositions of $R$ into indecomposable rings. See \cite[p.~100]{AndersonFuller}.

  The following proposition is an immediate consequence of Theorem~\ref{vhil}. 
  
 \begin{proposition}\label{3.8} Let $A$ be a right ideal of a ring $R$ such that $A$ has a serial factorization, but $R/A$ is not a uniserial right $R$-module.
 Then $A$ is a two-sided ideal,
 there exists a complete set of indecomposable central idempotents $\{x_1+A, \dots, x_n+A\}$ of $R/A$ and  $A=(A+(1-x_1)R)\dots(A+(1-x_n)R)$.\end{proposition}

From Proposition~\ref{3.8}, we get that:


\begin{proposition}\label{2.15} If a  right ideal $A$ of a ring $R$ has a serial factorization $A=A_1\dots A_n$ and $A$ has a set of generators of $n$ elements, then all the right ideals $A_1,\dots,A_n$ have a set of generators of $\le n+1$ elements.\end{proposition}

Assume that $A\subseteq B$ are two right ideals of a ring $R$. Suppose that $A$ has a serial factorization 
and that $B/A$ is a finitely generated right $R$-module. Then the right ideal $B$ can be generated by $A$ plus one further suitable 
element of~$B$ (Proposition~\ref{3.5}).

  \begin{corollary}\label{3.6} Let $R$ be a ring such that the right ideal $rR$ has a serial factorization for every non-zero $r\in R$. Then every finitely generated 
  non-zero right ideal $A$ of $R$ can be generated with two elements, the first of which can be any arbitrarily fixed non-zero element of $A$, and the 
  second must be suitably chosen. \end{corollary}
  
  \begin{Proof} Let $r\in A$ be any non-zero element. By Proposition~\ref{3.5}, $R/rR$ is a B\'ezout module.\end{Proof}

  In Commutative Algebra, the ideals $I$ of a ring $R$ with the property that every non-zero element $r\in I$ can be completed  to some two elements generating set of $I$ are some times called {\em one-and-a-half generated}. For instance, in a Dedekind domain, every ideal is one-and-a-half generated. Thus Corollary~\ref{3.6} says that if $R$ is a ring in which every principal right ideal has a serial factorization, then every finitely generated 
  right ideal of $R$ is one-and-a-half generated.
  
 
 
 
\bigskip

As we have already recalled in Example~\ref{2.4}(2), every non-zero ideal $A$ of a commutative Dedekind domain can be written as a product of powers of distinct maximal ideals $P_i$ in a unique way up to the order of the factors, $A=P_1^{t_1}\cdots P_k^{t_k}$. Every power $P_i^{t_i}$ is contained in a unique maximal ideal, which is $P_i$ itself. Moreover, every maximal ideal containing $A$ is one of these $n$ maximal ideals $P_i$. This fact holds, more generally, for our serial factorizations, as the following proposition shows.

\begin{proposition}\label{3.10} Let $A$ be a right ideal of a ring $R$ and suppose that $A$ has a serial factorization $A=A_1\dots A_n$. Then the following conditions hold:

{\rm (1)} Every right ideal $A_i$ is contained in a unique maximal right ideal $M_i$ of $R$.

{\rm (2)} The maximal right ideals $M_1,\dots,M_n$ are all distinct, and they are exactly the maximal right ideals of $R$ that contain $A$.

{\rm (3)} If $n\ge 2$, then the ideals $M_i$ are also maximal left ideals of $R$.\end{proposition}

\begin{Proof} Suppose $n=1$, so that $R/A$ is a uniserial right $R$-module. Then the serial factorization is the trivial factorization $A=A$, and (1) holds because $R/A$ has a unique maximal submodule. Now (2) is trivial.

Suppose $n\ge 2$. In this case, all the ideals $A$ and $A_i$ are
 two-sided (Theorem~\ref{vhil}), and $R/A\cong R/A_1\times\dots\times R/A_n$ canonically as rings, where the rings $R/A_i$ are right 
 chain rings. In this canonical isomorphism, the ideal $A_i/A$ of $R/A$ corresponds to the ideal $R/A_1\times\dots\times A_i/A_i\times\dots\times R/A_n$ and 
 the unique maximal right ideal of $R/A_1\times\dots\times R/A_n$ that contains it is $M'_i:=R/A_1\times\dots\times 
 J(R/A_i)\times\dots\times R/A_n$. The maximal right ideal of $R/A$ corresponding to $M_i'$ is $M_i/A$, where $M_i$ is the unique maximal right ideal of $R$ that contains $A_i$. The maximal 
 right ideals in the product of right chain rings $R/A_1\times\dots\times R/A_n$ are exactly the $n$ distinct maximal ideals $M'_1,
 \dots,M'_n$. In the canonically isomorphic ring $R/A$, we get that 
$M_1/A,\dots,M_n/A$ are exactly the $n$ distinct maximal ideals of $R/A$. Now (2) follows immediately. Finally, since the rings 
$R/A_i$ are right chain, the  ideals $M'_i$ are two-sided and 
also maximal left ideals, so that the ideals $M_i$ are also two-sided and maximal as left ideals.\end{Proof}

Recall that a commutative integral domain $S$ is {\em $h$-local} if each non-zero element of $S$ is contained in only a finite number of maximal ideals of $S$ and each non-zero prime ideal of $S$ is contained in only one maximal ideal of $S$ \cite[p.~27]{Matlis}. A {\em Pr\"ufer domain} is a commutative integral domain $R$ for which the localization $R_M$ of $R$ at $M$ is a valuation domain for every maximal ideal $M$ in $R$. Every commutative B\'ezout domain is a Pr\"ufer domain.

\begin{proposition}\label{h-local} {\rm (1)} If $R$ is a ring (possibly non-commutative), $r\in R$ and $rR$ has a serial factorization, then $r$ is contained in only a finite number of maximal right ideals of $R$ and each prime ideal $P$ of $R$ with $r\in P$ is contained in a unique maximal right ideal of $R$. 

{\rm (2)} A commutative integral domain $R$ has the property that every non-zero ideal has a serial factorization if and only if $R$ is an $h$-local Pr\"ufer domain.\end{proposition}

\begin{Proof} (1) The element $r$ is contained in only a finite number of maximal right ideals of $R$ by Proposition~\ref{3.10}(2).
If $r\in P$, then $rR=A_1\dots A_n\subseteq P$, so that $A_i\subseteq P$ for some $i$. But there is a unique maximal right ideal of $R$ that contains $A_i$, so that there is a unique maximal right ideal of $R$ that contains $P$.

(2) Let $I$ be a non-zero ideal of the Pr\"ufer $h$-local commutative domain $R$. Then $R/I$ is canonically isomorphic to $\oplus_M(R/I)_M$, where $M$ ranges in the set of all maximal ideals of $R$, equivalently in the set of all maximal ideals of $R$ that contain $I$ \cite[Theorem~22]{Matlis}. Let $M_1,\dots,M_n$ be the maximal ideals of $R$ that contain $I$. It follows that $R/I$ is canonically isomorphic to $\oplus_{i=1}^nR/(I_{M_i}\cap R)$. Let us show that $I=(I_{M_1}\cap R)\dots(I_{M_n}\cap R)$ is a serial factorization of $I$. The ideals $I_{M_i}\cap R$ are clearly proper. Since $I_{M_i}\cap R$ is contained only in the maximal ideal $M_i$, the ideals $I_{M_1}\cap R,\dots,I_{M_n}\cap R$ form a coindependent family. Thus $I=(I_{M_1}\cap R)\dots(I_{M_n}\cap R)$ by Lemma~\ref{1}. It remains to show that the modules $R/(I_{M_i}\cap R)$ are uniserial $R$-modules. They are cyclic modules over the localizations $R_{M_i}$, which are valuation domains, so that they are uniserial modules over $R_{M_i}$. In order to show that $R/(I_{M_i}\cap R)\cong (R/I)_{M_i}\cong R_{M_i}/I_{M_i}$ is uniserial as a module over $R$, it suffices to show that every cyclic $R$-submodule of $R_{M_i}/I_{M_i}$ is an $R_{M_i}$-submodule, that is, $xR_{M_i}+I_{M_i}\subseteq xR+I_{M_i}$ for every $x\in R_{M_i}$. Now, for every $t\in R\setminus M_i$, we have that $tR + (I_{M_i}\cap R) =R$, so that $tr+i=1$. Thus $x=trx+ix$, so $x+I_{M_i}=trx +I_{M_i}$ is divisible by $t$ in $xR+I_{M_i}/I_{M_i}$.

Conversely, let $R$ be a commutative integral domain in which every non-zero ideal has a serial factorization. Then $R$ is $h$-local by (1). In order to show that $R$ is Pr\"ufer, we must prove that the ideals of $R_M$ are linearly ordered by inclusion, for any fixed maximal ideal $M$. To this end, it suffices to prove that, for any non-zero ideal $J$ of $R_M$, the ideals of $R_M$ between $J$ and $R_M$ are linearly ordered by inclusion.
 Now the non-zero ideal $J$ of $R_M$ is extended from a non-zero ideal $I$ of $R$, that is, $J=I_M$. The ideal $I$ has  a serial factorization $I=I_1\dots I_n$, $R/I\cong R/I_1\oplus\dots\oplus R/I_n$, every  ideal $I_i$ is contained in a unique maximal 
 ideal $M_i$ of $R$, the maximal  ideals $M_1,\dots,M_n$ are all distinct, and they are exactly the maximal  ideals of $R$ that contain $I$ (Proposition~\ref{3.10}). Thus $I_i\cap (R\setminus M)\ne\emptyset$ for all the indices $i$ except for at most one index $i_0$. So $(R/I_i)_M=0$  for all the indices $i$ except for at most $i=i_0$. 
Hence $(R/I)_M\cong (R/I_{i_0})_M$ is the localization of the uniserial $R$-module $R/I_{i_0}$, hence it is a uniserial $R_M$-module. As $(R/I)_M\cong R_M/J$ is a uniserial $R_M$-module,
the ideals of $R_M$ between $J$ and $R_M$ are linearly ordered by set inclusion. 
\end{Proof}

Pr\"ufer $h$-local domains have received a lot of attention in the literature. See \cite{Matlis, Brandal, Olberding}.

\section{Right B\'ezout domains, right invariant elements and rigid elements}\label{4}

Let $R$ be a ring and $a$ be an element of $R$. Assume that the principal right ideal $aR$ has a serial factorization. By Theorem~\ref{vhil}, two cases can take place: (1) the trivial case where $R/aR$ is uniserial and the serial factorization is the trivial factorization $aR=aR$, or (2) $aR$ is a two-sided ideal of $R$. We will now consider this second case.
Clearly, $aR$ is a two-sided ideal of $R$ if and only if $Ra\subseteq aR$.

Now let $R$ be a (not necessarily commutative) integral domain. Recall that an element $a\in R$ is {\em right invariant} \cite[Section~0.9]{Cohn} if $Ra\subseteq aR$. {\em Left invariant} elements are
defined in a similar way, and an element $a$ is {\em invariant} if it is left and right invariant, that is, if $Ra= aR\ne 0$. The set Inv$(R)$ of all invariant elements of an integral domain $R$ is a multiplicatively closed subset of $R$ that contains all invertible elements of $R$. Notice that, in an integral domain, an element is right invertible if and only if it is left invertible. Also, if two elements of an equation
$x=yz$ between non-zero elements $x,y,z$ of the integral domain $R$
are invariant, then so is the third 
\cite[Lemma 0.9.1]{Cohn}.

An element $a$ of an integral domain $R$ is {\em rigid} if $a$ is non-zero, non-invertible, and for every $x,y,x'y'\in R$, $a=xy'=yx'$ 
implies $x = yu$ or $y = xu$ for some $u \in R$ (Equivalently,  for every $x,y,x'y'\in R$, $a=xy'=yx'$ 
implies $x' = uy'$ or $y' = ux'$ for some $u \in R$.)  See \cite[Section~0.7]{Cohn} and \cite[p.~2]{Zaf}.

\begin{Lemma}\label{vhip} Let $R$ be an integral domain and $a$ be a non-zero non-invertible element of $R$.

{\rm (1)} If $R/aR$ is uniserial, then $a$ is a rigid element of $R$.

{\rm (2)} If $R$ is a right B\'ezout domain, then $R/aR$ is uniserial if and only if $a$ is a rigid element.\end{Lemma}

\begin{Proof} (1) Suppose $R/aR$ uniserial and $a=xy'=yx'$. Then $aR$ is contained both in $xR$ and in $yR$. Since $R/aR$ is uniserial, it follows that either $xR\subseteq yR$ or $yR\subseteq xR$. Suppose for instance $xR\subseteq yR$. Then $x = yu$ 
 for some $u \in R$.
 
 (2) Let $R$ be a right B\'ezout domain and $a\in R$ be a rigid element. In order to show that $R/aR$ is uniserial, it suffices to show that any two cyclic submodules of $R/aR$ are comparable. Now a cyclic submodule of $R/aR$ is of the form $bR+aR/aR$. Since $R$ is right B\'ezout, every cyclic submodule of $R/aR$ is of the form $xR/aR$ with $xR\supseteq aR$, that is, of the form $xR/aR$ for some factorization $a=xy'$ of $a$.
Now it is easy to conclude that $R/aR$ is uniserial.\end{Proof}

In particular, in a right and left B\'ezout domain $R$, for a non-zero element $a$, $R/aR$ is uniserial if and only if $R/Ra$ is uniserial.

\bigskip

We say that two elements $a,b$ of an integral domain $R$ are {\em right associates} if there exists an invertible element $u\in R$ such that $a=bu$.

We will now generalize to the non-commutative setting the theory of semirigid GCD domains \cite{Zaf}. Recall that a commutative integral domain is a B\'ezout domain if and only if it is a Pr\"ufer GCD domain.

We say that an element $a$ of an integral domain $R$ is {\em semirigid} if $a$ is not rigid and $a$ has a factorization $a=a_1\dots a_n$ where:

(1) each $a_i$ is right invariant and rigid.

(2) $a_i$ and $a_j$ are right coprime (that is, $a_iR+a_jR=R$) for every $i\ne j$,

(3) $a_ia_j$ and $a_ja_i$ are right associates for every $i,j=1,2,\dots,n$.

\noindent
We will call such a factorization $a=a_1\dots a_n$ of a semirigid element $a\in R$ a {\em rigid factorization} of $a$.
Notice that the $n$ that appears in this factorization is necessarily $ \geq 2$. 
 (It would have been more precise to call our semirigid elements {\em right} semirigid elements, but we have preferred to simplify the terminology.)

\begin{proposition} \label{without} Let $a$ be a non-zero element of a right B\'ezout domain $R$. Then:

{\rm (1)} The principal right ideal $aR$ has a serial factorization if and only if either $a$ is invertible, or $a$ is rigid, or $a$ is semirigid. 

{\rm (2)} If $a=a_1\dots a_n$ is a rigid factorization of a semirigid element $a$, then $aR=a_1R\dots a_nR$ is a serial factorization of $aR$. 

{\rm (3)} Conversely, if $aR=A_1\dots A_n$ is a serial factorization of $aR$, where 
$n \geq 2$, then there exist $a_1,\dots,a_n\in R$ such that $A_i=a_iR$ for every $i$ and $a=a_1\dots a_n$ is a rigid factorization of $a$, so that $a$ is a semirigid element of $R$. \end{proposition}

\begin{Proof} Let $a$ be a non-zero non-invertible element of a right B\'ezout domain $R$. Suppose that the principal right ideal $aR$ has a serial factorization $aR=A_1\dots A_n$. If $n=1$, then $R/aR$ is uniserial, so that $a$ is rigid by Lemma~\ref{vhip}(2). Suppose $n\ge 2$. By Theorem~\ref{vhil}, all the right ideals $A_i$ and the right ideal $aR$ are two-sided, and $a$ is not rigid.
By Proposition~\ref{2.15}, each right ideal $A_i$ can be generated with two elements. Since $R$ is right B\'ezout, each right ideal $A_i$ is principal, $A_i=a_iR$ say, for a suitable non-zero element $a_i\in R$.  But $a_iR$ is a two-sided ideal. It follows that each $a_i$ is right invariant and $aR=a_1R\dots a_nR$. From this and the fact that all the elements $a_i$ are right invariant, we obtain that $a_1R\dots a_nR=a_1\dots a_nR$, because $a_1R\dots a_nR\subseteq a_1a_2R\dots a_rR\subseteq \dots\subseteq a_1\dots a_nR$. Thus $aR=a_1\dots a_nR$, which implies that $a=a_1\dots a_nu$ for some invertible element $u\in R$. Now $a_n$ rigid right invariant implies $a_nu$ rigid right invariant, because $Ra_nu\subseteq a_nRu=a_nR= a_nuR$. Since the ideals $a_iR$ are proper, the elements $a_i$ are not invertible. Changing the notation and denoting the element $a_nu$ by $a_n$, we get a factorization $a=a_1\dots a_n$ where each $a_i$ is right invariant and rigid. Condition (2) in the definition of rigid factorization follows trivially from the fact that the ideals $a_iR$ form a coindependent family of right ideals of $R$. Finally, $A_iA_j=A_jA_i$ implies that $a_iRa_jR=a_jRa_iR$, that is, $a_ia_jR=a_ja_iR$. It follows that there exists an invertible element $u_{ij}\in R$ such that $a_ia_j=a_ja_iu_{ij}$. This proves one implication in statement (1) and statement (3).

To conclude the proof, it suffices to show that (2) holds. Suppose that $a$ is semirigid and $a=a_1\dots a_n$ is a rigid factorization of $a$. In order to show that $aR=a_1R\dots a_nR$ is a serial factorization of $aR$, the only non-trivial things to be checked is that the right ideals $a_iR$ form a coindependent family. 
 This follows from Lemma~\ref{2.2}(1).\end{Proof}

If $R$ is a right B\'ezout domain, $a$ is a non-zero element of $R$ and $u$ is an invertible element of $R$, then:

\noindent (1) $a$ is rigid if and only if $ua$ is rigid (because $R/aR$ is isomorphic to $R/uaR$ via left multiplication by $u$).

\noindent (2) If $a$ is right invariant, then $ua$ is right invariant. (Because, for every $r\in R$, $rua=u(u^{-1}ru)a=uas$ for some $s\in R$.)

\noindent (3) If $a$ is semirigid, then $ua$ is semirigid and $ua$ is right associate to $a$ (by Theorem~\ref{3.3}).

Let us pass to the uniqueness of a rigid factorization of a semirigid element. From Theorem~\ref{vhil}, we immediately get that:

\begin{theorem}\label{4.3} Let $R$ be a right B\'ezout domain and $a\in R$ be a semirigid element. Let $a=a_1\dots a_n=b_1\dots b_m$ be two rigid factorizations of $A$. Then $n=m$ and there exists a unique permutation $\sigma$ of $\{1,\dots,n\}$ such that $a_i$ and $b_{\sigma(i)}$ are right associates for every $i=1,\dots,n$. \end{theorem}

Theorem~\ref{4.3} extends \cite[Theorem~2]{Zaf} to the case of non-commutative rings.

\begin{proposition}\label{lattices} Let $a=a_1\dots a_n$ be a rigid factorization of a semirigid element $a$ of a right B\'ezout domain $R$. Then any right invariant left divisor $b$ of $a$, that is, any right invariant element $b\in R$ for which there exists $c\in R$ with $a=bc$, factorizes as $b=b_1\dots b_n$, where each $b_i$ is a right invariant left divisor of $a_i$. \end{proposition} 

\begin{Proof} Apply Theorem \ref{Tuesday} to the right ideals $A=aR$ and $B=bR$ of $R$, obtaining that the serial factorization of $bR$ is $bR=(bR+a_1R)\dots(bR+a_nR)$ (where we are supposed to omit the factors $bR+a_iR$ equal to $R$). Taking as $b_i$ a (suitable) generator of the two-sided ideal $bR+a_iR$, 
we get the desired factorization $b=b_1\dots b_n$ of $b$.
\end{Proof}


For every right invariant element $a$ of $R$, we can consider the lattice $L_{ri}(aR,R)$ whose elements are all the principal right ideals $bR$ of $R$, where $b$ is any right invariant left divisor of $a$. According to Proposition~\ref{lattices}, if $a=a_1\dots a_n$ is a rigid factorization of a semirigid element $a$ of a right B\'ezout domain $R$, then the 
lattice $L_{ri}(aR,R)$ is canonically isomorphic to the direct product $\prod_{i=1}^nL_{ri}(a_iR,R)$ of the $n$ chains $L_{ri}(a_iR,R)$.

\medskip

We conclude our paper with a variation of Proposition \ref{lattices}:

\begin{corollary} Let $a=a_1\dots a_n$ be a rigid factorization of a semirigid element $a$ of a right B\'ezout domain $R$. Then any non-invertible right invariant left divisor $b$ of $a$  
is either rigid or has a rigid factorization  $b=b_{i_1}\dots b_{i_t}$ for suitable integers $t\ge 1,\ 1\le i_1<i_2<\dots<i_t\le n$, where each $b_{i_j}$ is a left divisor of $a_{i_j}$.\end{corollary}

\begin{Proof} As in the proof of Proposition \ref{lattices}, apply Theorem \ref{Tuesday} to the right ideals $A=aR$ and $B=bR$ of $R$, obtaining that the serial factorization of $bR$ is $bR=(bR+a_1R)\dots(bR+a_nR)$. Let $1\le i_1<i_2<\dots<i_t\le n$ be the indices with $bR+a_iR$ a proper right ideal, so that $bR=(bR+a_{i_1}R)\dots(bR+a_{i_t}R)$. 
If $b$ is not rigid, then $t \geq 2$ and we can apply Proposition \ref{without}(3) to get to the 
desired result. 
 \end{Proof}

We are grateful to the referee for a very careful reading of our manuscript and to the Centro de Matem\'atica of the University of Porto for its hospitality when we were writing part of this paper.

 {}

\end{document}